\documentclass[12pt]{article}

\usepackage{amsmath}
\usepackage{amsthm}
\usepackage{amssymb}
\usepackage{amscd}
\usepackage{a4wide}
\newtheorem{dfn}{Definition}[section]
\newtheorem{thm}[dfn]{Theorem}

\newtheorem{cor}[dfn]{Corollary}
\newtheorem{lemma}[dfn]{Lemma}

\newtheorem{remark}[dfn]{\it Remark}

\newcommand\GL{\mathrm{GL}}
\newcommand\Symp{\mathrm{Sp}}
\newcommand\ep{\varepsilon}
\newcommand\trans{{}^t\!}
\newcommand\vectx{\boldsymbol{x}}
\newcommand\vecty{\boldsymbol{y}}
\newcommand\vecta{\boldsymbol{a}}
\newcommand\vectb{\boldsymbol{b}}
\newcommand\Nat{\mathbb{N}}
\newcommand\Comp{\mathbb{C}}
\newcommand\Real{\mathbb{R}}
\renewcommand\tilde{\widetilde}

\title{%
A bialternant formula for odd symplectic characters \\
and its application
}

\author{Soichi Okada%
\footnote{
Graduate School of Mathematics, Nagoya University, 
Furo-cho, Chikusa-ku, Nagoya 464-8602, Japan
{\tt okada@math.nagoya-u.ac.jp}
}
}

\date{
\it Dedicated to Professor Toshio Oshima on the occasion of his 70th birthday
}

\begin{document}

\maketitle

\begin{abstract}
We present a bialternant formula for odd symplectic characters, 
which are the characters of indecomposable modules of odd symplectic groups 
introduced by R.~Proctor. 
As an application, we give a linear algebraic proof to an odd symplectic character identity 
due to R.~P.~Brent, C.~Krattenthaler and S.~O.~Warnaar.
\end{abstract}

\section{%
Introduction
}

The irreducible characters of classical (semisimple or reductive) groups are expressed 
as a ratio of two Vandermonde-type determinants 
(with modifications needed for even special orthogonal groups). 
This bialternant formula is a restatement of Weyl's character formula, 
and Weyl's denominator formula gives a factorization of the determinant in the denominator.
These determinant formulas paly a key role in the representation theory of classical groups 
and the combinatorial applications. 
One can derive many classical group character identities such as 
the Jacobi--Trudi identity, the Giambelli identity and the Cauchy identity 
(see \cite{NSSY}, \cite{O98} and \cite{O06} for example).
The aim of this paper is to give a bialternant formula for odd symplectic groups, 
which are neither semisimple nor reductive.

A \emph{partition} is a weakly decreasing sequence 
$\lambda = (\lambda_1, \lambda_2, \dots)$ of nonnegative integers 
such that $|\lambda| = \sum_i \lambda_i$ is finite.
The \emph{length} $l(\lambda)$ of a partition $\lambda$ is defined to 
be the number of positive parts of $\lambda$.
We denote by $\varnothing$ the empty partition $(0,0,\dots)$.

Recall Weyl's character and denominator formulas for symplectic groups $\Symp_{2n}(\Comp)$.
For a partition $\lambda$ of length $\le n$, 
we denote by $\Symp_{2n}(\lambda ; x_1, \dots, x_n)$ the irreducible character 
of $\Symp_{2n}(\Comp)$ corresponding to $\lambda$.
Then Weyl's character formula is written in bialternant form as
\begin{equation}
\label{eq:Sp-char}
\Symp_{2n}(\lambda ; x_1, \dots, x_n)
 =
\frac{ \det \Big( x_i^{\lambda_j+n-j+1} - x_i^{-(\lambda_j+n-j+1)} \Big)_{1 \le i, j \le n} }
     { \det \Big( x_i^{n-j+1} - x_i^{-(n-j+1)} \Big)_{1 \le i, j \le n} },
\end{equation}
and Weyl's denominator formula reads
\begin{multline}
\label{eq:Sp-den}
\det \Big( x_i^{n-j+1} - x_i^{-(n-j+1)} \Big)_{1 \le i, j \le n}
\\
 =
\prod_{i=1}^n
 \big( x_i - x_i^{-1} \big)
\prod_{1 \le i < j \le n}
 \big( x_i^{1/2} x_j^{1/2} - x_i^{-1/2} x_j^{-1/2} \big)
 \big( x_i^{1/2} x_j^{-1/2} - x_i^{1/2} x_j^{-1/2} \big).
\end{multline}

The symplectic group $\Symp_{2n}(\Comp)$ is the subgroup of $\GL_{2n}(\Comp)$ 
consisting of elements preserving a non-degenerate skew-symmetric bilinear form.
Non-degenerate skew-symmetric bilinear forms exist only in even dimension.
Proctor \cite{P} introduced ``symplectic groups on odd dimensional vector spaces'' 
and developed the representation/character theory of them.
Proctor's odd symplectic group $\Symp_{2n+1}(\Comp)$ is defined to be the subgroup 
of $\GL_{2n+1}(\Comp)$ consisting of elements preserving 
a skew-symmetric bilinear form of maximal possible rank.
He obtained an indecomposable $\Symp_{2n+1}(\Comp)$-module $V_\lambda$ 
for each partition $\lambda$ of length $\le n+1$, 
and gave several formulas for the character $\Symp_{2n+1}(\lambda ; x_1, \dots, x_n ; z)$ 
of $V_\lambda$, which we call the \emph{odd symplectic character}.
See Definition~\ref{def:osp-char} for a representation-free definition of 
odd symplectic characters.

The main result of this article is the following bialternant formula 
for the odd symplectic characters.

\begin{thm}
\label{thm:Weyl}
For a partition $\lambda$ of length $\le n+1$, 
we define $A_\lambda = A_\lambda(x_1, \dots, x_n ; z)$ 
to be the $(n+1) \times (n+1)$ matrix with $(i,j)$ entry given by
$$
\begin{cases}
 \big( x_i^{\lambda_j+n-j+2} - x_i^{-(\lambda_j+n-j+2)} \big)
- z^{-1} \big( x_i^{\lambda_j+n-j+1} - x_i^{-(\lambda_j+n-j+1)} \big)
 &\text{if $1 \le i \le n$,} \\
z^{\lambda_j+n-j+1}
 &\text{if $i = n+1$.}
\end{cases}
$$
Then we have
\begin{multline}
\label{eq:osp-den}
\det A_{\varnothing}
\\
 =
\prod_{i=1}^n \big( x_i - x_i^{-1} \big)
\prod_{1 \le i < j \le n+1}
 \big( x_i^{1/2} x_j^{1/2} - x_i^{-1/2} x_j^{-1/2} \big)
 \big( x_i^{1/2} x_j^{-1/2} - x_i^{-1/2} x_j^{1/2} \big),
\end{multline}
where $x_{n+1} = z$. 
And we have
\begin{equation}
\label{eq:osp-char}
\Symp_{2n+1}(\lambda ;x_1, \dots, x_n ; z)
 =
\frac{ \det A_{\lambda} }
     { \det A_{\varnothing} }.
\end{equation}
\end{thm}

If we put $z = 1$ in this theorem, we can recover Proctor's bialternant formula 
for $\Symp_{2n+1}(\lambda; x_1, \dots, x_n ; 1)$ (\cite[Theorem~2.2]{P}).
We note that the factors of the denominator $\det A_{\varnothing}$ correspond to the set of vectors
$$
\{ 2 \ep_i :  1 \le i \le n \}
\cup
\{ \ep_i \pm \ep_j : 1 \le i < j \le n+1 \},
$$
which sits between the positive system of type $C_n$ and that of type $C_{n+1}$, 
where $(\ep_1, \dots, \ep_{n+1})$ is the standard orthonormal basis of $\Real^{n+1}$.

As an application of the bialternant formula for odd symplectic characters, 
we give a linear algebraic proof of the Brent--Krattenthaler--Warnaar identity \cite{Kr2},
which were found in their study of discrete Mehta-type integrals \cite{BKW}
and proved by using a combinatorial interpretation of odd symplectic characters.
Their identity is an ``odd symplectic'' analogue of the identities given in 
\cite[Theorem~2.2]{O98}.

\begin{thm}
\label{thm:BKW}
\textup{(}Brent--Krattenthaler--Warnaar \cite{Kr2}\textup{)}
Let $m$ and $n$ be positive integers with $m \le n$, and $r$ a nonnegative integer.
Then we have
\begin{multline}
\label{eq:BKW}
\sum_\lambda
 z^{-r}
 \,
 \Symp_{2m+1}(\lambda ; x_1, \dots, x_m ; z)
 \,
 \Symp_{2n+1}((r^{n-m}) \cup \lambda ; y_1, \dots, y_n ; z)
 \\
=
\Symp_{2(m+n+1)}( (r^{m+n+1}) ; x_1, \dots, x_m, y_1, \dots, y_n, z),
\end{multline}
where $\lambda$ runs over all partitions with $l(\lambda) \le m+1$, 
and $(r^{n-m}) \cup \lambda$ and $(r^{n+m+1})$ stand for the partitions 
$(\underbrace{r, \dots, r}_{n-m}, \lambda_1, \dots, \lambda_m)$ 
and $(\underbrace{r, \dots, r}_{m+n+1})$ respectively.
\end{thm}

This paper is organized as follows.
In Section~2, we define odd symplectic characters in terms of generating function 
and give a proof of our main theorem (Theorem~\ref{thm:Weyl}).
Section~3 is devoted to the proof of Brent--Krattenthaler--Warnaar identity (Theorem~\ref{thm:BKW}).
 % Introduction
\section{%
Bialternant formula
}

In this section we define odd symplectic characters 
in terms of generating function (Cauchy-type formula) 
and give a proof of Theorem~\ref{thm:Weyl}.

Recall that the \emph{Schur function} $s_\lambda(x_1, \dots, x_n)$ 
corresponding to a partition $\lambda$ of length $\le n$ is given by 
\begin{equation}
\label{eq:GL-char}
s_\lambda(x_1, \dots, x_n)
 =
\frac{ \det \Big( x_i^{\lambda_j+n-j} \Big)_{1 \le i, j \le n} }
     { \det \Big( x_i^{n-j} \Big)_{1 \le i, j \le n} }.
\end{equation}
The Schur functions are the characters of irreducible polynomial representations 
of the general linear group $\GL_n(\Comp)$ and the above definition (\ref{eq:GL-char}) 
is nothing but Weyl's character formula. 
And Weyl's denominator formula for $\GL_n(\Comp)$ is the evaluation of the Vandermonde determinant:
\begin{equation}
\label{eq:GL-den}
\det \Big( x_i^{n-j} \Big)_{1 \le i, j \le n}
 =
\prod_{1 \le i < j \le n} \big( x_i - x_j \big).
\end{equation}
In this article we define odd symplectic characters in terms of generating function 
with respect to Schur functions.

\begin{dfn}
\label{def:osp-char}
\textup{(}See \cite[Proposition~3.1]{P}\textup{)}
Let $x_1, \dots, x_n$ and $z$ be indeterminates.
We define odd symplectic characters $\Symp_{2n+1}(\lambda;  x_1, \dots, x_n ; z)$ 
for partitions $\lambda$ of length $\le n+1$ by the relation
\begin{multline}
\label{eq:osp-Cauchy}
\sum_\lambda
 \Symp_{2n+1}(\lambda; x_1, \dots, x_n ; z)
 s_\lambda(u_1, \dots, u_{n+1})
\\
=
\frac{ \prod_{1 \le i < j \le n+1} (1 - u_i u_j) }
     { \prod_{i=1}^n \prod_{j=1}^{n+1} (1 - x_i u_j) (1 - x_i^{-1} u_j)
       \prod_{j=1}^{n+1} (1 - z u_j) },
\end{multline}
where $\lambda$ runs over all partitions of length $\le n+1$.
\end{dfn}

Since the Schur functions are linearly independent, 
the odd symplectic characters are uniquely determined by (\ref{eq:osp-Cauchy}).
It follows from (\ref{eq:osp-Cauchy}) that $\Symp_{2n+1}(\lambda; x_1, \dots, x_n ; z)$ 
is a Laurent polynomial in $x_1, \dots, x_n$ and a polynomial in $z$.
Note that the odd symplectic characters $\Symp_{2n+1}(\lambda ; x_1, \dots, x_n ; z)$ 
are obtained from symplectic universal character $\Symp(\lambda;\vectx)$ 
by specializing $\vectx = (x_1, \dots, x_n, z, x_1^{-1}, \dots, x_n^{-1}, 0, 0, \dots)$ 
(see \cite{Ki}, \cite{KT} and \cite{O06}).

Recall the Cauchy--Binet formula and the Cauchy determinant, which play a fundamental role 
in this paper.
Given an $M \times N$ matrix $X = (x_{i,j})_{1 \le i \le M, 1 \le j \le N}$ 
and a subset $I \subset [N] = \{ 1, \dots, N \}$ of column indices, 
we denote by $X(I)$ the submatrix of $X$ obtained by picking up columns indexed by $I$.
Then the Cauchy--Binet formula is stated as follows:

\begin{lemma}[Cauchy--Binet formula]
\label{lem:CauchyBinet}
For two $M \times N$ matrices $X$ and $Y$, we have
$$
\sum_I \det X(I) \det Y(I)
 =
\det \left( \trans X Y \right),
$$
where $I$ runs over all $M$-element subsets of $[N]$.
\end{lemma}

The following determinant evaluations are known as the Cauchy determinants.

\begin{lemma}[Cauchy determinants]
\label{lem:Cauchy}
For indeterminates $x_1, \dots, x_n$ and $y_1, \dots, y_n$, we have
\begin{gather}
\label{eq:Cauchy1}
\det \left(
 \frac{ 1 }
      { x_i - y_j }
\right)_{1 \le i, j \le n}
 =
\frac{ (-1)^{n(n-1)/2} \prod_{1 \le i < j \le n} (x_i - x_j) (y_i - y_j) }
     { \prod_{i=1}^n \prod_{j=1}^n (x_i - y_j) },
\\
\label{eq:Cauchy2}
\det \left(
 \frac{ 1 }
      { 1 - x_i y_j }
\right)_{1 \le i, j \le n}
 =
\frac{ \prod_{1 \le i < j \le n} (x_i - x_j) (y_i - y_j) }
     { \prod_{i=1}^n \prod_{j=1}^n (1 - x_i y_j) }.
\end{gather}
\end{lemma}

For a partition $\lambda$ of length $\le n+1$, we put
$$
I_{n+1}(\lambda) = \{ \lambda_1 + n, \lambda_2 + n-1, \dots, \lambda_n + 1, \lambda_{n+1} \}.
$$
Then the correspondence $\lambda \mapsto I_{n+1}(\lambda)$ 
gives a bijection from the set of partitions of length $\le n+1$ to the set of $(n+1)$-element subsets 
of $\Nat$, the set of nonnegative integers.

By using these lemmas, we prove Theorem~\ref{thm:Weyl}.

\begin{proof}[Proof of Theorem~\ref{thm:Weyl}]
First we prove (\ref{eq:osp-den}).
By adding the $i$th column multiplied by $z^{-1}$ to the $(i-1)$st column for $i=n, n-1, \dots, 2$, 
and then by multiplying the last row by $z-z^{-1}$, 
we obtain
$$
\det A_{\varnothing}
 =
\frac{ 1 }
     { z - z^{-1} }
\det \Big( x_i^{n-j} - x_i^{-(n-j)} \Big)_{1 \le i, j \le n+1},
$$
where we put $x_{n+1} = z$.
By pulling out $x_i-x_i^{-1}$ from the $i$th row 
and then by using elementary column operations, we see that
$$
\det \Big( x_i^{n-j} - x_i^{-(n-j)} \Big)_{1 \le i, j \le n+1}
 =
\prod_{i=1}^{n+1} \big( x_i - x_i^{-1} \big)
\cdot
\det \Big( \big( x_i + x_i^{-1} \big)^{n-j+1} \Big)_{1 \le i, j \le n+1}.
$$
Now we can obtain (\ref{eq:osp-den}) by applying the Vandermonde determinant and 
using
$$
\big( x_i + x_i^{-1} \big) - \big( x_j + x_j^{-1} \big)
 =
\big( x_i^{1/2} x_j^{1/2} - x_i^{-1/2} x_j^{-1/2} \big)
\big( x_i^{1/2} x_j^{-1/2} - x_i^{-1/2} x_j^{1/2} \big).
$$

Next we prove (\ref{eq:osp-char}).
Since the odd symplectic characters $\Symp_{2n+1}(\lambda ; x_1, \dots, x_n ; z)$ are 
characterized by (\ref{eq:osp-Cauchy}), 
it is enough to show that the ratios $\det A_{\lambda} / \det A_{\varnothing}$ satisfy 
\begin{equation}
\label{eq:osp-Cauchy2}
\sum_\lambda
 \frac{ \det A_{\lambda} }
      { \det A_{\varnothing} }
 s_\lambda(u_1, \dots, u_{n+1})
 =
\frac{ \prod_{1 \le i < j \le n+1} (1 - u_i u_j) }
     { \prod_{i=1}^n \prod_{j=1}^{n+1} (1 - x_i u_j) (1 - x_i^{-1} u_j)
       \prod_{j=1}^{n+1} (1 - z u_j) }.
\end{equation}

Let $X = (x_{i,k})_{1 \le i \le n+1, k \ge 0}$ 
and $Y = (y_{i,k})_{1 \le i \le n+1, k \ge 0}$ be the $(n+1) \times \infty$ matrices 
whose $(i,k)$ entries are given by
\begin{align*}
x_{i,k}
 &=
\begin{cases}
x_i^{k+1} - x_i^{-k-1} - z^{-1} (x_i^k - x_i^{-k}) &\text{if $1 \le i \le n$,} \\
z^k &\text{if $i=n+1$,}
\end{cases}
\\
y_{i,k}
 &=
u_i^k.
\end{align*}
Then we have
$$
\frac{ \det A_{\lambda} }
     { \det A_{\varnothing} }
=
\frac{ \det X(I_{n+1}(\lambda)) }
     { \det X(I_{n+1}(\varnothing)) },
\quad
s_\lambda(u_1, \dots, u_{n+1})
 =
\frac{ Y(I_{n+1}(\lambda)) }
     { Y(I_{n+1}(\emptyset)) }.
$$
Straightforward computations show that the $(i,j)$ entry of $X \trans Y$ is equal to
$$
\begin{cases}
\dfrac{ (x_i - x_i^{-1}) ( 1 - z^{-1} u_j ) }
      { ( 1 - x_i u_j) ( 1 - x_i^{-1} u_j) }
 &\text{if $1 \le i \le n$,}
\\
\dfrac{ 1 }{ 1 - z u_j }
 &\text{if $i=n+1$.}
\end{cases}
$$
By applying the Cauchy--Binet formula (Lemma~\ref{lem:CauchyBinet}) 
and then by pulling out common factors, 
we obtain
\begin{multline*}
\sum_\lambda
 \det X(I_{n+1}(\lambda)) \det Y(I_{n+1}(\lambda))
\\
 =
\prod_{i=1}^n \big( x_i - x_i^{-1} \big)
\prod_{j=1}^{n+1} \big( 1 - z^{-1} u_j \big)
\cdot
\det \left(
\frac{ 1 }
     { ( 1 - x_i u_j ) ( 1 - x_i^{-1} u_j ) }
\right)_{1 \le i, j \le n+1},
\end{multline*}
where $\lambda$ runs over all partitions of length $\le n+1$ and $x_{n+1} = z$.
By substituting $x_i + x_i^{-1}$ and $u_j + u_j^{-1}$ for $x_i$ and $y_j$ respectively 
in the Cauchy determinant (\ref{eq:Cauchy1}), we have
\begin{multline*}
\det \left( \frac{ 1 }{ ( 1 - x_i u_j ) ( 1 - x_i^{-1} u_j ) } \right)_{1 \le i, j \le n+1}
\\
 =
\prod_{1 \le i < j \le n+1}
 \big( x_i^{1/2} x_j^{1/2} - x_i^{-1/2} x_j^{-1/2} \big)
 \big( x_i^{1/2} x_j^{-1/2} - x_i^{-1/2} x_j^{1/2} \big)
\\
\times
\frac{ \prod_{1 \le i < j \le n+1} (u_i - u_j)(1 - u_i u_j) }
     { \prod_{i=1}^{n+1} \prod_{j=1}^{n+1} ( 1 - x_i u_j )( 1 - x_i^{-1} u_j ) }.
\end{multline*}
Therefore, combining these computations together with (\ref{eq:osp-den}) 
and the Vandermonde determinant (\ref{eq:GL-den}), we obtain (\ref{eq:osp-Cauchy2}) 
and complete the proof of (\ref{eq:osp-char}).
\end{proof}

By putting $z = 1$ in Theorem~\ref{thm:Weyl}, we can recover Proctor's bialternant formula.

\begin{cor}
\textup{(}Proctor \cite[Theorem~2.2 and Proposition~7.1]{P}\textup{)}
For a partition of length $\le n+1$, let $B_\lambda = B_\lambda(x_1, \dots, x_n)$ 
be the $(n+1) \times (n+1)$ matrix whose $(i,j)$ entry is given by
$$
\begin{cases}
 x_i^{\lambda_j+n-j+3/2} + x_i^{-(\lambda_j+n^j+3/2)}
 &\text{if $1 \le i \le n$,}
\\
 1 &\text{if $i=n+1$.}
\end{cases}
$$
Then we have
\begin{multline}
\label{eq:osp-den2}
\det B_{\varnothing}
 =
\prod_{i=1}^n \big( x_i^{1/2} - x_i^{-1/2} \big) \big( x_i - x_i^{-1} \big)
\\
\times
\prod_{1 \le i < j \le n}
  \big( x_i^{1/2} x_j^{1/2} - x_i^{-1/2} x_j^{-1/2} \big)
  \big( x_i^{1/2} x_j^{-1/2} - x_i^{-1/2} x_j^{1/2} \big),
\end{multline}
and
\begin{equation}
\label{eq:osp-char2}
\Symp_{2n+1}(\lambda;x_1, \dots, x_n;1)
 =
\frac{ \det B_{\lambda} }
     { \det B_{\varnothing} }.
\end{equation}
\end{cor}

\begin{proof}[Proof]
The proof of (\ref{eq:osp-den2}) is similar to that of (\ref{eq:osp-den}), so we omit it.
By substituting $z = 1$, the $(i,j)$ entry of $A_\lambda$ becomes
\begin{multline*}
\big( x_i^{\lambda_j+n-j+2} - x_i^{-(\lambda_j+n-j+2)} \big)
 -
\big( x_i^{\lambda_j+n-j+1} - x_i^{-(\lambda_j+n-j+1)} \big)
\\
 =
\big( x_i^{1/2} - x_i^{-1/2} \big)
\big( x_i^{\lambda_j+n-j+3/2} + x_i^{-(\lambda_j+n-j+3/2)} \big)
\end{multline*}
for $1 \le i \le n$.
Hence by cancelling the common factors $\prod_{i=1}^n (x_i^{1/2} - x_i^{-1/2})$ 
we obtain (\ref{eq:osp-char2}) from (\ref{eq:osp-char}).
\end{proof}

\begin{remark}
By specializing $z=-1$ in \textup{(\ref{eq:osp-char})}, 
we obtain Krattenthaler's bideterminantal formula \cite[(3.5)]{Kr1} 
for $\Symp_{2n+1}(\lambda ; x_1, \dots, x_n ; -1)$.
\end{remark}

Appealing to Weyl's denominator formula for type $D_{n+1}$, 
we can derive a product formula for the principal specialization of odd symplectic characters 
from Proctor's bialternant formula (\ref{eq:osp-char2}).
Let $\Delta^+(D_{n+1})$ be the positive system of type $D_{n+1}$ given by
$$
\Delta^+(D_{n+1})
 =
\{ \ep_i \pm \ep_j : 1 \le i < j \le n+1 \},
$$
where $(\ep_1, \dots, \ep_{n+1})$ is the standard orthonormal basis of $\Real^{n+1}$ 
with respect to the inner product $\langle \ , \ \rangle$.

\begin{cor}
For a partition $\lambda$ of length $\le n+1$, we have
\begin{equation}
\label{eq:osp-q}
\Symp_{2n+1}(\lambda ; q^n, q^{n-1}, \dots, q ; 1)
 =
\prod_{\alpha \in \Delta^+(D_{n+1})}
 \frac{ \big[ \langle \lambda + \rho, \alpha \rangle \big]_q }
      { \big[ \langle \rho, \alpha \rangle \big]_q },
\end{equation}
where $\rho = (n+1/2, n-1/2, \dots, 3/2, 1/2)$ and $[m]_q = (q^{m/2} - q^{-m/2})/(q^{1/2} - q^{-1/2})$.
\end{cor}

\begin{proof}[Proof]
Weyl's denominator formula for type $D_{n+1}$ reads
\begin{multline*}
\det \Big( x_i^{n-j+1} + x_i^{-(n-j+1)} \Big)_{1 \le i, j \le n+1}
\\
 =
2 
\prod_{1 \le i < j \le n+1}
  \big( x_i^{1/2} x_j^{1/2} - x_i^{-1/2} x_j^{-1/2} \big)
  \big( x_i^{1/2} x_j^{-1/2} - x_i^{-1/2} x_j^{1/2} \big).
\end{multline*}
By specializing $x_i = q^{\lambda_i+n-i+3/2}$ for $1 \le i \le n+1$, we obtain 
$$
\det B_\lambda( q^n, q^{n-1}, \dots, q )
 =
\prod_{\alpha \in \Delta^+(D_{n+1})} \big[ \langle \lambda+\rho, \alpha \rangle \big]_q.
$$
Hence (\ref{eq:osp-q}) follows from (\ref{eq:osp-char2}).
\end{proof}
 % Bialternant formula
\section{%
Application to Brent--Krattenthaler--Warnaar's identity
}

In this section, we use the bialternant formula (Theorem~\ref{thm:Weyl}) 
to prove Brent--Krattenthaler--Warnaar's identity (Theorem~\ref{thm:BKW}).
The idea of our proof is the same as \cite[Theorem~2.2]{O98}.
Before the proof, we prepare two lemmas.
The first one enables us to reduce the proof of Theorem~\ref{thm:BKW} 
to the case where $m=n$.

\begin{lemma}
\label{lem:reduction}
\begin{enumerate}
\item[\textup{(a)}]
Let $\lambda$ be a partition with length $\le n+1$ and $\lambda_1 \le r$.
Then
$$
(x_1 \cdots x_n)^r \, \Symp_{2n+1}(\lambda ; x_1, \dots, x_n ; z)
$$
is a polynomial in $x_1, \dots, x_n$, 
and we have
\begin{multline*}
\big[
(x_1 \cdots x_n)^r \, \Symp_{2n+1}(\lambda; x_1, \dots, x_n ; z)
\big]
\big|_{x_1 = 0}
\\
 =
\begin{cases}
(x_2 \cdots x_n)^r \, \Symp_{2n-1}( (\lambda_2, \dots, \lambda_n) ; x_2, \dots, x_n ; z)
 &\text{if $\lambda_1 = r$,} \\
0 &\text{otherwise},
\end{cases}
\end{multline*}
where $F|_{x_1=0}$ means that we substitute $x_1 = 0$ in $F$.
\item[\textup{(b)}]
Let $\lambda$ be a partition with length $\le n$ and $\lambda_1 \le r$.
Then
$$
(x_1 \cdots x_n)^r \, \Symp_{2n}(\lambda ; x_1, \dots, x_n)
$$
is a polynomial in $x_1, \dots, x_n$, 
and we have
\begin{multline*}
\big[
(x_1 \cdots x_n)^r \, \Symp_{2n}(\lambda; x_1, \dots, x_n)
\big]
\big|_{x_1 = 0}
\\
 =
\begin{cases}
(x_2 \cdots x_n)^r \, \Symp_{2(n-1)}( (\lambda_2, \dots, \lambda_n) ; x_2, \dots, x_n)
 &\text{if $\lambda_1 = r$,} \\
0 &\text{otherwise}.
\end{cases}
\end{multline*}
\end{enumerate}
\end{lemma}

\begin{proof}[Proof]
(a)
It follows from Theorem~\ref{thm:Weyl} that
$$
(x_1 \cdots x_n)^r \, \Symp_{2n+1}(\lambda ; x_1, \dots, x_n ; z)
=
\frac{ (x_1 \cdots x_n)^{r+n+1} \det A_{\lambda} }
     { (x_1 \cdots x_n)^{n+1} \det A_{\varnothing} }.
$$
Let $\tilde{A}_\lambda = \big( \tilde{a}_{i,j} \big)_{1 \le i, j \le n+1}$ 
be the $(n+1) \times (n+1)$ matrix whose $(i,j)$ entry is given by
$$
\tilde{a}_{i,j}
 =
x_i^{r+n+1}
\Big(
 \big( x_i^{\lambda_j+n+2-j} - x_i^{-\lambda_j-n-2+j} \big)
 - z^{-1} \big( x_i^{\lambda_j+n+1-j} - x_i^{-\lambda_j-n-1+j} \big)
\Big)
$$
for $1 \le i \le n$ and $\tilde{a}_{n+1,j} = z^{\lambda_j+n+1-j}$.
Then we have
$$
(x_1 \cdots x_n)^{r+n+1} \det {A_\lambda}
 =
\det \tilde{A}_\lambda.
$$
Since $\lambda_1 \le r$ by assumption, we have $\lambda_j \le r$ for all $1 \le j \le n+1$.
If $1 \le i \le n$ and $1 \le j \le n+1$, we have
$$
\tilde{a}_{i,j}
 =
x_i^{\lambda_j+r+2n-j+3} - x_i^{r-\lambda_j+j-1}
 - z^{-1} x_i^{\lambda_j+r+2n-j+2} + z^{-1} x_i^{r-\lambda_j+j}
$$
with
$$
\lambda_j+r+2n-j+3 > \lambda_j+r+2n-j+2 \ge r-\lambda_j+j > r-\lambda_j+j-1 \ge 0.
$$
Hence we see that $(x_1 \cdots x_n)^r \, \Symp_{2n+1}(\lambda; x_1, \dots, x_n ; z)$ 
is a polynomial in $x_1, \dots, x_n$.
Since $r - \lambda_j+j-1 = 0$ if and only if $j=1$ and $\lambda_1 = r$, we have
$$
\tilde{a}_{1,j} \big|_{x_1=0}
 =
\begin{cases}
 -1 &\text{if $\lambda_1 = r$ and $j=1$,} \\
 0 &\text{otherwise.}
\end{cases}
$$
Hence we have
$$
\big( \det \tilde{A}_\lambda \big) \big|_{x_1=0}
 =
\begin{cases}
 (-1) \cdot \det \big( \tilde{a}_{i+1,j+1} \big)_{1 \le i, j \le n}
 &\text{if $\lambda_1 = r$,}
\\
 0
 &\text{if $\lambda_1 < r$.}
\end{cases}
$$
If $\lambda_1 = r$, then we have
$$
\det \big( \tilde{a}_{i+1,j+1} \big)_{1 \le i, j \le n}
 =
(x_2 \cdots x_n)^{r+n+1} \det A_{(\lambda_2, \dots, \lambda_{n+1})} (x_2, \dots, x_n;z).
$$
Therefore we have
\begin{align*}
&
\big[
 (x_1 \cdots x_n)^r \, \Symp_{2n+1}(\lambda ; x_1, \dots, x_n ; z)
\big]
\big|_{x_1=0}
\\
&\quad
 =
\frac{ (x_2 \cdots x_n)^{r+n+1} \, \det A_{(\lambda_2, \dots, \lambda_{n+1})}(x_2, \dots, x_n ; z) }
     { (x_2 \cdots x_n)^{n+1} \, \det A_{\varnothing}(x_2, \dots, x_n ; z) }
\\
&\quad
 =
(x_2 \cdots x_n)^r \, \Symp_{2n-1}( (\lambda_2, \dots, \lambda_n) ; x_2, \dots, x_n ; z ).
\end{align*}
This completes the proof of (a).

The proof of (b) is similar to (a). (See \cite[Lemma~5.3]{O98}.)
\end{proof}

The following lemma is a key to our proof of Theorem~\ref{thm:BKW}.

\begin{lemma}
\label{lem:key}
Let $p(x,y,z,a,b)$ be the rational function in $x$, $y$, $z$, $a$ and $b$ 
given by
\begin{multline}
\label{eq:p}
p(x,y,z,a,b)
\\
 =
\frac{ (1-xz)(1-yz) }{ 1 - xy }
-
a \cdot \frac{ (x-z)(1-yz) }{ x-y }
+
b \cdot \frac{ (1-xz)(y-z) }{ x-y }
-
ab \cdot \frac{ (x-z)(y-z) }{ 1-xy}.
\end{multline}
Let $\vectx = (x_1, \dots, x_n)$, $\vecty = (y_1, \dots, y_n)$, 
$\vecta = (a_1, \dots, a_n)$, $\vectb = (b_1, \dots, b_n)$, $z$ and $c$ be indeterminates.
We define an $(n+1) \times (n+1)$ matrix $C = C(\vectx,\vecty,z;\vecta,\vectb,c) 
= \big( C_{i,j} \big)_{1 \le i, j \le n+1}$ by putting
$$
C_{ij}
=
\begin{cases}
p(x_i,y_j,z,a_i,b_j) &\text{if $1 \le i, \, j \le n$,} \\
1 - a_i &\text{if $i=n+1$ and $1 \le j \le n$,} \\
1 - b_j &\text{if $1 \le i \le n$ and $j=n+1$,} \\
\dfrac{ 1 - c }{ 1 - z^2 } &\text{if $i = j = n+1$,}
\end{cases}
$$
and a $(2n+1) \times (2n+1)$ matrix $V = V(\vectx,\vecty,z ; \vecta,\vectb,c) 
= \big( V_{i,j} \big)_{1 \le i, j \le 2n+1}$ by putting
$$
V_{i,j} = x_i^{j-1} - a_i x_i^{2n+1-j},
\quad
V_{n+i,j} = y_i^{j-1} - b_i y_i^{2n+1-j},
\quad
V_{2n+1,j} = z^{j-1} - c z^{2n+1-j}
$$
for $1 \le i \le n$ and $1 \le j \le 2n+1$.
Then we have
\begin{equation}
\label{eq:key}
\det C
 =
\frac{ (-1)^n }
     { (1 - z^2) \prod_{i=1}^n \prod_{j=1}^n (x_i - y_j)(1 - x_i y_j) } 
\det V.
\end{equation}
\end{lemma}

\begin{proof}[Proof]
Since the both sides of (\ref{eq:key}) are polynomials in $a_1, \dots, a_n, b_1, \dots, b_n$ 
with degree at most one in each variable, 
they are linear combinations of $a^I b^J = \prod_{i \in I} a_i \prod_{j \in J} b_j$ 
for subsets $I$ and $J \subset [n]$.
We denote by $L(I,J)$ and $R(I,J)$ the coefficients of $a^I b^J$ 
on the left and right hand sides respectively, 
and prove $L(I,J) = R(I,J)$.

Fix two subsets $I$ and $J$ of $[n]$.
Let $\sigma = \sigma_{I,J}$ be the involutive ring automorphism of the Laurent polynomial ring defined by
$$
\sigma(x_i)
 =
\begin{cases}
x_i^{-1} &\text{if $i \in I$,} \\
x_i      &\text{if $i \not\in I$,}
\end{cases}
\quad
\sigma(y_j)
 =
\begin{cases}
y_j^{-1} &\text{if $j \in J$,} \\
y_j      &\text{if $j \not\in J$.}
\end{cases}
$$
We compute $\sigma(L(I,J))$ and $\sigma(R(I,J))$ explicitly.

It follows from (\ref{eq:p}) that $L(I,J)$ is equal to the determinant of the matrix 
$C' = \big( C'_{i,j} \big)_{1 \le i, j \le n+1}$ whose $(i,j)$ entry is given by
$$
C'_{ij}
=
\begin{cases}
 - \dfrac{ (x_i-z)(y_j-z) }{ 1-x_iy_j } &\text{if $i \in I$ and $j \in J$,} \\
 - \dfrac{ (x_i-z)(1-y_jz) }{ x_i-y_j } &\text{if $i \in I$ and $j \in [n] \setminus J$,} \\
 \dfrac{ (1-x_iz)(y_j-z) }{ x_i-y_j } &\text{if $i \in [n] \setminus I$ and $j \in J$, } \\
 \dfrac{ (1-x_iz)(1-y_jz) }{ 1 - x_iy_j } &\text{if $i \in [n] \setminus I$ and $j \in [n] \setminus J$,} \\
 -1 &\text{if $i \in I$ and $j = n+1$,} \\
 1 &\text{if $i \in [n] \setminus I$ and $j=n+1$,} \\
 -1 &\text{if $i=n+1$ and $j \in J$,} \\
 1 &\text{if $i=n+1$ and $j \in [n] \setminus J$,} \\
 \dfrac{ 1-c }{ 1 - z^2 } &\text{if $i=j=n+1$.}
\end{cases}
$$
Hence, for $1 \le i, j \le n$, we have
$$
\sigma(C'_{i,j})
 =
\begin{cases}
 \dfrac{ (1-x_iz)(1-y_jz) }{ 1-x_iy_j } &\text{if $i \in I$ and $j \in J$,} \\
 - \dfrac{ (1-x_iz)(1-y_jz) }{ 1-x_iy_j } &\text{if $i \in I$ and $j \in [n] \setminus J$,} \\
 - \dfrac{ (1-x_iz)(1-y_jz) }{ 1-x_iy_j } &\text{if $i \in [n] \setminus I$ and $j \in J$, } \\
 \dfrac{ (1-x_iz)(1-y_jz) }{ 1-x_iy_j } &\text{if $i \in [n] \setminus I$ and $j \in [n] \setminus J$.}
\end{cases}
$$
By pulling out common factors from rows and columns of $\sigma(C')$, we obtain
$$
\det \sigma(C')
 =
(-1)^{\# I + \# J}
\prod_{i=1}^n (1 - x_i z) \prod_{j=1}^n (1 - y_j z)
\cdot \det C'',
$$
where the $(i,j)$ entry of the matrix $C'' = \big( C''_{i,j} \big)_{1 \le i, j \le n+1}$ is given by
$$
C''_{i,j}
 =
\begin{cases}
 \dfrac{ 1 }{ 1 - x_i y_j } &\text{if $1 \le i \le j$ and $1 \le j \le n$,} \\
 \dfrac{ 1 }{ 1 - x_i z} &\text{if $1 \le i \le n$ and $j=n+1$,} \\
 \dfrac{ 1 }{ 1 - y_j z} &\text{if $i=n+1$ and $1 \le j \le n$,} \\
 \dfrac{ 1 - c }{ 1 - z^2 } &\text{if $i=j=n+1$.}
\end{cases}
$$
Here we can use the Cauchy determinant (\ref{eq:Cauchy2}) to evaluate 
the coefficients of $c^0$ and $c^1$ in $\det C''$:
\begin{align*}
[c^0] \det C''
 &=
\frac{ \prod_{1 \le i < j \le n} (x_i - x_j) (y_i - y_j)
       \prod_{i=1}^n (x_i - z) (y_i - z) }
     { \prod_{i,j=1}^n (1 - x_i y_j)
       \prod_{i=1}^n (1 - x_i z) (1 - y_j z)
       (1 - z^2) },
\\
[c^1] \det C''
 &=
-
\frac{ 1 }{ 1 - z^2 }
\cdot
\frac{ \prod_{1 \le i < j \le n} (x_i - x_j) (y_i - y_j) }
     { \prod_{i=1}^n \prod_{j=1}^n ( 1 - x_i y_j) }.
\end{align*}
Therefore we have
\begin{multline}
\label{eq:LHS}
\sigma(L(I,J))
 =
(-1)^{\# I + \# J}
\frac{ \prod_{1 \le i < j \le n} (x_i - x_j) (y_i - y_j) }
     { (1 - z^2) \prod_{i,j=1}^n (1 - x_i y_j) }
\\
\times
\left[
 \prod_{i=1}^n (x_i - z) (y_i - z)
 - c \prod_{i=1}^n (1 - x_i z) (1 - y_i z)
\right].
\end{multline}

Next we compute $\sigma( R(I,J) )$.
We see that
$$
R(I,J)
 = 
\frac{ (-1)^n }
     { (1 - z^2) \prod_{i=1}^n \prod_{j=1}^n (x_i - y_j)(1 - x_i y_j) } 
\det V',
$$
where the entries of $V' = \big( V'_{i,j} \big)_{1 \le i, j \le 2n+1}$ are given by
\begin{align*}
V'_{i,j}
 &=
\begin{cases}
 - x_i^{2n+1-j} &\text{if $i \in I$,} \\
 x_i^{j-1} &\text{if $i \in [n] \setminus I$,}
\end{cases}
\\
V'_{n+i,j}
 &=
\begin{cases}
 - y_i^{2n+1-j} &\text{if $j \in J$,} \\
 y_i^{j-1} &\text{if $j \in [n] \setminus J$,}
\end{cases}
\\
V'_{2n+1,j}
 &= z^{j-1} - c z^{2n+1-j}.
\end{align*}
Since we have
$$
\sigma \big( (x_i - y_j)(1 - x_i y_j) \big)
 =
\begin{cases}
x_i^{-2} y_j^{-2} (x_i - y_j) (1 - x_i y_j) &\text{if $i \in I$ and $j \in J$,} \\
x_i^{-2} (x_i - y_j) (1 - x_i y_j) &\text{if $i \in I$ and $j \in [n] \setminus J$,} \\
y_j^{-2} (x_i - y_j) (1 - x_i y_j) &\text{if $i \in [n] \setminus I$ and $j \in J$,} \\
(x_i - y_j) (1 - x_i y_j) &\text{if $i \in [n] \setminus I$ and $j \in [n] \setminus J$,}
\end{cases}
$$
we obtain
\begin{align*}
&
\sigma(R(I,J))
\\
&=
(-1)^n
\frac{ \prod_{i \in I} x_i^{2n} \prod_{j \in J} y_j^{2n} }
     { (1 - z^2) \prod_{i=1}^n \prod_{j=1}^n (x_i - y_j)(1 - x_i y_j) }
\det \sigma(V')
\\
&=
(-1)^{n + \# I + \# J}
\frac{ 1 }
     { (1 - z^2) \prod_{i=1}^n \prod_{j=1}^n (x_i - y_j)(1 - x_i y_j) }
\det V'',
\end{align*}
where $V''$ is the $(2n+1) \times (2n+1)$ matrix whose $(i,j)$ entry $V''_{i,j}$ is given by
$$
V''_{i,j} = x_i^{j-1},
\quad
V''_{n+i,j} = y_i^{j-1},
\quad
V''_{2n+1,j} = z^{j-1} - c z^{2n-j+1}.
$$
Here we can use the Vandermonde determinant (\ref{eq:GL-den}) to evaluate the coefficients 
of $c^0$ and $c^1$ in $\det V''$:
\begin{align*}
[c^0] \det V''
 &=
\prod_{1 \le i < j \le n} (x_j - x_i) (y_j - y_i)
\prod_{i,j = 1}^n (y_j - x_i)
\prod_{i=1}^n (z - x_i)(z - y_i),
\\
[c^1] \det V''
 &=
(-1) z^{2n}
\prod_{1 \le i < j \le n} (x_j - x_i)(y_j - y_i)
\prod_{i,j=1}^n (y_j - x_i)
\\
&\quad
\times
\prod_{i=1}^n (z^{-1} - x_i)(z^{-1} - y_j).
\end{align*}
Therefore we have
\begin{multline}
\label{eq:RHS}
\sigma(R(I,J))
 =
(-1)^{\# I + \# J}
\frac{ \prod_{1 \le i < j \le n} (x_i - x_j) (y_i - y_j) }
     { (1 - z^2) \prod_{i,j=1}^n (1 - x_i y_j) }
\\
\times
\left[
 \prod_{i=1}^n (x_i - z) (y_i - z)
 - c \prod_{i=1}^n (1 - x_i z) (1 - y_i z)
\right].
\end{multline}

By (\ref{eq:LHS}) and (\ref{eq:RHS}), we have $\sigma(L(I,J)) = \sigma(R(I,J))$ 
and this completes the proof of Lemma~\ref{lem:key}.
\end{proof}

Now we are ready to give a proof of Theorem~\ref{thm:BKW}.

\begin{proof}[Proof of Theorem~\ref{thm:BKW}]
First we prove the case where $m=n$.
We put $M = r+n$.
We apply the Cauchy--Binet formula (Lemma~\ref{lem:CauchyBinet}) to the matrices 
$X = \big( X_{i,k} \big)_{1 \le i \le n+1, 0 \le k \le M}$ and 
$Y = \big( Y_{i,k} \big)_{1 \le i \le n+1, 0 \le k \le M}$ whose entries are given by
\begin{align*}
X_{i,k}
 &=
\begin{cases}
x_i^{k+1} - x_i^{-k-1} - z^{-1} (x_i^k - x_i^{-k}) &\text{if $1 \le i \le n$,} \\
z^k &\text{if $i=n+1$,}
\end{cases}
\\
Y_{i,k}
 &=
\begin{cases}
y_i^{k+1} - y_i^{-k-1} - z^{-1} (y_i^k - y_i^{-k}) &\text{if $1 \le i \le n$,} \\
z^k &\text{if $i=n+1$.}
\end{cases}
\end{align*}
It follows from a direct but lengthy computation that the $(i,j)$ entry of $X \trans Y$ is given by
$$
\begin{cases}
x_i^{-(M+1)} y_j^{-(M+1)} p(x_i,y_j,z^{-1},x_i^{2M+2},y_j^{2M+2})
 &\text{if $1 \le i \le n$ and $1 \le j \le n$},
\\
- x_i^{-(M+1)} z^M \big( 1 - x_i^{2M+2} \big)
 &\text{if $1 \le i \le n$ and $j = n+1$,}
\\
- y_j^{-(M+1)} z^M \big( 1 - y_j^{2M+2} \big)
 &\text{if $i = n+1$ and $1 \le j \le n$,}
\\
\dfrac{ 1 - z^{2(M+1)} }
      { 1 - z^2 }
 &\text{if $i = j = n+1$.}
\end{cases}
$$
By pulling out factors $x_i^{-(M+1)}$ from the $i$th row, $y_j^{-(M+1)}$ from the $j$th column 
and $-z^M$ from the last row and column, 
we see that
$$
\det X \trans Y
 =
\prod_{i=1}^n x_i^{-M-1} \prod_{j=1}^n y_j^{-M-1} z^{2M}
\cdot
\det C(\vectx, \vecty, z^{-1} ; \vectx^{2M+2}, \vecty^{2M+2}, z^{-2M-2}),
$$
where $C$ is the matrix introduced in Lemma~\ref{lem:key}, 
$\vectx^{2M+2} = (x_1^{2M+2}, \dots, x_n^{2M+2})$
and $\vecty^{2M+2} = (y_1^{2M+2}, \dots, y_n^{2M+2})$.
Hence, by using Theorem~\ref{thm:Weyl} and applying the Cauchy--Binet formula (Lemma~\ref{lem:CauchyBinet}), 
we have
\begin{align*}
&
\sum_\lambda
z^{-r}
\, 
\Symp_{2n+1}(\lambda ; x_1, \dots, x_n ; z)
\, 
\Symp_{2n+1}(\lambda ; y_1, \dots, y_n ; z)
\\
&=
\sum_{l(\lambda) \le n+1}
z^{-r}
\frac{ \det X(I_{n+1}(\lambda)) }{ \det X(I_{n+1}(\varnothing)) }
\cdot
\frac{ \det Y(I_{n+1}(\lambda)) }{ \det Y(I_{n+1}(\varnothing)) }
\\
&=
z^{-r} \frac{ 1 }{ \det X(I_{n+1}(\varnothing)) \det Y(I_{n+1}(\varnothing)) }
\det (X \trans Y)
\\
&=
z^{-r} \frac{ 1 }{ \det X(I_{n+1}(\varnothing)) \det Y(I_{n+1}(\varnothing)) }
\\
&\quad\quad\times
\prod_{i=1}^n x_i^{-M-1} y_i^{-M-1} z^{2M}
\cdot
\det C(\vectx, \vecty, z^{-1} ; \vectx^{2M+2}, \vecty^{2M+2}, z^{-2M-2}),
\end{align*}
where $\lambda$ runs over all partitions of length $\le n+1$ satisfying $\lambda_1 \le r$.
Now we can use Lemma~\ref{lem:key} to obtain
\begin{align*}
&
\sum_\lambda
z^{-r}
\,
\Symp_{2n}(\lambda ; x_1, \dots, x_n ; z)
\,
\Symp_{2n}(\lambda ; y_1, \dots, y_n ; z)
\\
&=
z^{-r}
\frac{ 1 }{ \det X(I_{n+1}(\varnothing)) \det Y(I_{n+1}(\varnothing)) }
\cdot
\frac{ (-1)^n 
       \prod_{i=1}^n x_i^{-M-1} y_i^{-M-1} z^{2M} }
     { (1 - z^{-2}) \prod_{i=1}^n \prod_{j=1}^n (x_i - y_j)(1 - x_i y_j) } 
\\
&\quad\quad\times
\det V(\vectx,\vecty,z^{-1} ; \vectx^{2M+2}, \vecty^{2M+2}, z^{-2M-2}).
\end{align*}
Here $\det X(I_{n+1}(\varnothing)) = (-1)^{n(n+1)/2} \det A_{\varnothing}$ 
and $\det Y(I_{n+1}(\varnothing))$ are evaluated 
by using (\ref{eq:osp-den}) and the Vandermonde determinant (\ref{eq:GL-den}) respectively.
Also we have
\begin{multline*}
\det V(\vectx,\vecty,z^{-1} ; \vectx^{2M+2}, \vecty^{2M+2}, z^{-2M-2})
\\
=
\prod_{i=1}^n (-1) x_i^{r+2n+1}
\prod_{j=1}^n (-1) y_j^{r+2n+1}
\cdot
z^{-r-2n-1}
\\
\times
\det \left( t_i^{r+2n+2-j} - t_i^{-(r+2n+2-j)} \right)_{1 \le i, j \le 2n+1},
\end{multline*}
where $(t_1, \dots, t_{2n+1}) = (x_1, \dots, x_n, y_1, \dots, y_n, z)$.
Therefore, by using (\ref{eq:Sp-char}) and (\ref{eq:Sp-den}), we conclude that
\begin{multline*}
\sum_\lambda
z^{-r}
\,
\Symp_{2n+1}(\lambda ; x_1, \dots, x_n ; z)
\,
\Symp_{2n+1}(\lambda ; y_1, \dots, y_n ; z)
 \\=
\Symp_{4n+2}((r^{2n+1}) ; x_1, \dots, x_n, y_1, \dots, y_n, z).
\end{multline*}
This complete the proof of Theorem~\ref{thm:BKW} in the case $m=n$.

Lastly we prove the general case by the downward induction on $m$.
We assume (\ref{eq:BKW}).
By multiply the both sides of (\ref{eq:BKW}) by $(x_1 \dots x_m y_1 \dots y_n z)^r$, 
and then by substituting $x_1 = 0$, it follows from Lemma~\ref{lem:reduction} that
\begin{multline*}
\sum_\lambda
 z^{-r}
 \, 
 (x_2 \dots x_m)^r \, \Symp_{2m-1}( (\lambda_2, \dots, \lambda_{m+1}) ; x_2, \dots, x_m ; z)
\\
\times
 (y_1 \dots y_n z)^r \, \Symp_{2n+1}( (r^{n-m}) \cup \lambda ; y_1, \dots, y_n ; z)
\\
=
(x_2 \dots x_m y_1 \dots y_n z)^r \, \Symp_{2(m+n)}( (r^{m+n}) ; x_1, \dots, x_m, y_1, \dots, y_n, z),
\end{multline*}
where the summation is taken over all partitions $\lambda$ of length $\le m+1$ such that $\lambda_1 = r$.
Since $(r^{n-m}) \cup \lambda = (r^{n-m+1}) \cup (\lambda_2, \dots, \lambda_{n+1})$ 
for such a partition, we obtain
\begin{multline*}
\sum_\mu
 z^{-r}
 \,
 \Symp_{2m-1}( \mu ; x_2, \dots, x_m ; z )
 \,
 \Symp_{2n+1}( (r^{n-m+1}) \cup \mu ; y_1, \dots, y_n ; z)
\\
=
\Symp_{2(m+n)}( (r^{m+n}) ; x_2, \dots, x_m, y_1, \dots, y_n, z),
\end{multline*}
where $\mu$ runs over all partitions of length $m$ satisfying $\mu_1 \le r$.
This is the desired identity, which is Equation (\ref{eq:BKW}) with $m$ replaced by $m-1$.
\end{proof}
 % Application to BKW identity
\subsection*{Acknowledgements}
This work was partially supported by Grant-in-Aid for Scientific Research 
No.~24340003 and No.~18K03208.
The author gratefully acknowledges the support and hospitality 
of the Galileo Galilei Institute for Theoretical Physics (Italy) and 
the National Institute for Mathematical Sciences (Korea), 
where part of this work was carried out.
He thanks Christian Krattenthaler for sharing unpublished work 
on an odd symplectic character identity.
 % Acknowledgements
%% Please do NOT use ``\bysame'' command in your bibliography list.

 % References

\end{document}